# On Separation of Variables

C. P. Viazminsky
Department of Physics
University of Aleppo
Syria

**Abstract.**
The necessary and sufficient conditions for a function $F(x_1,....,x_n)$ to be totally or partially separable are derived. It is shown that a function is totally separable if and only if each component of the gradient vector of $Ln|F|$ depends only on the corresponding variable. The conditions of separability are expressed neatly in terms of the matrix $F_{,ij} - F_{,i}F_{,j}$ which has to be diagonal if the function is to be totally separable, and has to assume a diagonal block form in order that the function is partially separable. The conditions of separability are also given without using derivatives. For polynomials, the conditions of separability are shown to hold if and only if the product of the first column and the first row of the coefficients matrix is equal to the matrix itself. This promotes an easy computational scheme for checking, and actually carrying out, variable separation.

## 1. Introduction

Separation of variables in a function is a common practice in special types of ordinary and partial differential equations [1,2]. Separation of variables plays an essential role is statistical distributions, where independence of random variables is equivalent to the possibility of factorizing the common distribution function to products of distributions each of which depends on one random variable [3]. The process of factorization of polynomials was studied in [4], where a separation scheme was designed. The latter scheme however, requires function's derivatives for its implementation.

Before getting involved in the technicalities it may be helpful to give some basic general facts underlying the subject of this work. A function $F(x_1,....,x_n)$, defined on a region $D \subseteq \Re^n$, is said to be separable with respect to the variable $x_1$ (or in the variable $x_1$) if it can be written as a product of two functions, defined on D, with one function depending only on $x_1$, whereas the other is independent of $x_1$. The function $F(x_1,....,x_n)$ is said totally separable ( or just separable) if it was separable in all its variables. When separable, the function $F(x_1,....,x_n)$ can be factorized into a product of $n$ functions $f_i(x_i)$ that depend on distinct variables. A function $F(x_1,....,x_n)$ is called partially separable if it was separable with some, but not all, of its variables. The process in which the function $F$ is factorized to products that depend on distinct variables is called variable separation. As mentioned earlier, variable separation could be partial or total. It is clear that the factors $f_i(x_i)$ appearing in variable separation of a function $F$ are unique up to arbitrary constant multipliers $\alpha_i$, such that $\prod_{i=1}^{n} \alpha_i = 1$.

We finally point out to the close link between additive separability and multiplicative separability. A function $F(x_1,....,x_n)$ is additively separable, or separable with respect to addition, if it can be written in the form
$$F(x_1,....,x_n) = g_1(x_1) + .... + g_n(x_n)$$
It clear that: $F$ is separable with respect to addition if and only if $\exp(F)$ is separable with respect to multiplication. Or equivalently, $F$ is separable with respect to multiplication if and only if $Ln|F|$ is separable with respect to addition.

# 1. Separation of One Variable

We assume throughout this section that the function $F: D \subseteq \Re^n \to \Re$ is differentiable with respect to $x_1$ and its partial derivative $\partial F / \partial x_1 \equiv F_{,1}$ is continuous on D.

**Theorem 1**. The function $F(x_1,....,x_n)$ is separable in its first variable if and only if
$$F_{,1}/F \tag{2.1}$$
depends only on $x_1$.

Proof. ($\Rightarrow$) If $F$ is separable with respect to $x_1$, it can be written in the form
$$F = f_1(x_1)R_1(x_2,...,x_n), \tag{2.2}$$
and hence
$$F_{,1}/F = f_{1,1}/f_1 . \tag{2.3}$$
As the right hand-side of eq(2.3) is a function of $x_1$ alone, so must be its left hand-side.

($\Leftarrow$) Assume that $F_{,1}/F$ is a function $X_1$ of $x_1$ alone:
$$F_{,1}/F = X_1(x_1). \tag{2.4}$$
Integrating both sides of this equation with respect to $x_1$ yields
$$Ln|F/R_1(x_2,....,x_n)| = \int X_1(x_1)dx_1 \tag{2.5}$$
where $R_1$ is an arbitrary function which does not depends on $x_1$. Solving (2.5) for $F$ yields
$$F = f_1(x_1)R_1(x_2,....,x_n), \tag{2.6}$$
where
$$|f_1(x_1)| = \exp(\int X_1(x_1)dx_1). \tag{2.7}$$
If $R_1(x_2,....,x_n)$ is left arbitrary, then (2.6) yields all functions that satisfy (2.4). For our objective function $F(x_1,....,x_n)$, the function $R_1$ is determined by
$$R_1(x_2,....,x_n) = F(x_1,....,x_n)/f_1(x_1). \tag{2.8}$$

Note that we could have solved eq(2.4), using the standard method of integrating a linear partial differential equation, to obtain the same general solution (2.6).

# 3. Conditions of Total Separability

**Theorem 2**. Assume that the function $F(x_1,....,x_n)$ is of class $C^1$ on its domain of definition $D \subseteq \Re^n$. In order that $F$ is totally separable, it is necessary and sufficient that

$$F_{,i}/F \qquad (3.1)$$

depends on $x_i$ alone, for $i = 1,....,n$.

Proof. ($\Rightarrow$) Assume that $F$ is separable in all its variables, i.e. it can be factorized in the form:

$$F = f_1(x_1)....f_n(x_n). \qquad (3.2)$$

It is evident that

$$F_{,i}/F = f_{i,i}/f_i \qquad (3.3)$$

depends on $x_i$ alone. This is true for $i = 1,....,n$.

($\Leftarrow$) Assume that $F_{,i}/F$ is dependent on $x_i$ only, where $i = 1,....,n$, and set

$$F_{,i}/F = X_i(x_i). \qquad (3.4)$$

It is easy to check that the simultaneous system of linear partial differential equations (3.4) has the following integrals

$$\begin{aligned}F(x_1,....,x_n) &= f_1(x_1)R_1(x_2,....,x_n) \\ &= f_2(x_2)R_2(x_1,x_3,....,x_n) \\ &\vdots \\ &= f_n(x_n)R_n(x_1,....,x_n).\end{aligned} \qquad (3.5)$$

From the last equations it follows that, for each value of $i$, the ratio

$$F/f_i(x_i) = R_i(x_1,...,x_{i-1},x_{i+1},...,x_n) \qquad (3.6)$$

is independent of $x_i$. Multiplying the function $F/f_1$, which is independent of $x_1$, by $1/f_2$, yields $F/f_1 f_2$ which is certainly still independent of $x_1$. Moreover, and if $F/f_1 f_2$ depends on $x_2$, then its product by $f_1(x_1)$, namely, $f_1(F/f_1 f_2) = F/f_2$ must depends on $x_2$, which is in contradiction with (3.5). Thus $F/f_1 f_2$ is independent of both $x_1$ and $x_2$. The proof proceeds by induction to conclude that $F/f_1 f_2 .... f_n = const$, because this ratio is independent of all the variables $x_1,....,x_n$. Thus

$$F = const.f_1.f_2. .... .f_n, \qquad (3.7)$$

and the proof is complete.¦

**Remarks**:
1. Had we left the constant multiplier in (3.7) arbitrary we would have obtained all functions that satisfy the system of differential equations (3.4).
2. Since separability of a function $F(x_1,....,x_n)$ with respect to its first $(n-1)$ variables implies its total separabilty, a second version of theorem 2 is:

**Theorem 2a**. A function $F(x_1,....,x_n)$ of class $C^1(D)$ is totally separable if and only if $F_{,i}/F$ is a function of $x_i$ alone, for any set of $n-1$ chosen variables.

3. Theorem (2) may be deduced as a direct generalization of theorem 1. Indeed, the function $F = f_1 R_1$, which is already separable with respect to $x_1$, is also separable with respect to $x_2$, i.e. it is has the form $F = f_1 f_2 R_{12}(x_3,....,x_n)$, if and only if $F_{,2}/F$ depends on $x_2$ alone. The proof proceeds to conclude that $F$ is separable with respect to $x_1,....,x_{n-1}$, and thus having the form :

$$F = f_1....f_{n-1}R_{12...n-1}(x_n) \qquad (3.8)$$

if and only if the conditions (3.1) are valid for $i = 1,...,n-1$. The last form however shows that $F$ is separable in all its variables. We may thus take $f_n \equiv R_{12...n-1}(x_n)$, so that the constant in (3.7) is 1.

4. The conditions of separability (3.4) may be expressed compactly as

$$\nabla Ln|F| = (X_1(x_1),...., X_n(x_n)), \qquad (3.9)$$

which is equivalent to

$$dLn|F| = \sum_{i=1}^{n} X_i(x_i)dx_i , \qquad (3.10)$$

which in turn is equivalent to

$$|F| = c.\exp(\sum_{i=1}^{n} \int X_i(x_i)dx_i)$$
$$\text{i.e. } F = c.\prod_{i=1}^{n} f_i(x_i). \qquad (3.11)$$

The content of this remark constitutes another proof of theorem 2 which may be rephrased in the following way.

**Theorem 2b.** A function $F$ of class $C^1(D)$ is totally separable if and only if each component of $\nabla Ln|F|$ depends only on the corresponding variable.

If $F$ is separable, the discussion above together with equation (3.4), promote a scheme for its separation by the formula

$$F = const.\prod_{i=1}^{n} \exp \int \frac{F_{,i}}{F} dx_i \qquad (3.12)$$

This scheme however, is not practical, since it requires performing a number of order $3n$ operations of partial differentiation to calculate $F_{,i}$, divisions to calculate $F/F_{,i}$ and integrations. We shall find much easier separation schemes in the forthcoming sections.

**Theorem 3.** A function $F(x_1,....,x_n)$ of class $C^2(D)$ is totally separable if and only if the following conditions hold

$$FF_{,ij} - F_{,i}F_{,j} = 0 \qquad (i \neq j). \qquad (3.13)$$

Proof. By theorem (2) $F$ is separable if and only if, for each $i$, $(F_{,i}/F)$ depends on $x_i$ alone. These conditions are equivalent to the following

$$(F_{,i}/F)_{,j} = 0 \qquad (i \neq j).$$

The theorem follows on noting that

$$(\frac{F_{,i}}{F})_{,j} = \frac{F_{,ij}F - F_{,i}F_{,j}}{F^2} \qquad (3.14)$$

On the account of the symmetry of equations (3.13) with respect to its indices, the number of independent conditions that are to be satisfied by a function $F$, in order to be separable, reduces to $\frac{1}{2}n(n-1)$ conditions corresponding to

$$FF_{,ij} - F_{,i}F_{,j} = 0 \quad (i < j). \tag{3.15}$$

The latter conditions have to be proved all valid if the function $F$ is to be separable. For $n = 2$, we obtain Scott's theorem [1,2] which states that a function $F(x,y)$ is separable if and only if

$$FF_{,xy} - F_{,x}F_{,y} = 0.$$

## 4. Separation With Respect to Sets of Variables

Let $I$ be a proper subset of the set $\{1,2,\ldots n\}$, and $J = \{1,\ldots,n\} - I$. The variables $x_i \, (i \in I)$ will be denoted collectively by $x_I$, so that a function $G$ that depends only on the variables $x_i \, (i \in I)$ will be denoted by $G(x_I)$. Similarly, $H(x_J)$ will denote a function which depends only on $x_j \, (j \in J)$.

**Theorem 4**. In order that a function $F(x_1,\ldots,x_n)$ of class $C^1(D)$ is separable in the set of variables $x_I$, and hence capable of assuming the form,

$$F = G(x_I).H(x_J), \tag{4.1}$$

it is necessary and sufficient that the functions

$$\frac{F_{,i}}{F} \equiv g_i(x_I) \quad (\forall i \in I). \tag{4.2}$$

depend only on $x_I$.

Proof. First, it is clear that if $F$ is separable in the form (4.1), then conditions (4.2) are satisfied, and hence are necessary. In proving the converse we suppose that the given function $F$ satisfies the conditions (4.2). Let us look temporarily on the variables $x_j \, (j \in J)$ in $F$ as parameters, and thus denote the differential of $Ln|F|$, under this assumption, by $d_I Ln|F|$. Now

$$d_I Ln|F| = \sum_{i \in I} \frac{F_{,i}}{F} dx_i$$

By hypothesis, the right hand side does not involve the variables $x_J$, and hence the left hand side, when $x_J$ are considered parameters, also does not involve $x_J$. Thus the right hand-side is a total differential, which we designate by $dg(x_I)$, because the left hand-side is so. Therefore

$$dLn|F| = dg(x_I) + \sum_{j \in J} \frac{F_{,j}}{F} dx_j .$$

Since $dLn|F| - dg$ is a total differential, the sum on the right hand-side must be a total differential, which is designated by $dh(x_J)$. Therefore $Ln|F| = g(x_I) + h(x_J)$, and $F = G(x_I)H(x_J)$, where $g(x_I) = Ln|G(x_I)|$, and $h(x_J) = Ln|H(x_J)|$. ¦

**Corollary 1**. It is clear that a function $F$ is separable in $x_I$ if and only if it is separable in $x_J$, and hence, iff

$$\frac{F_{,j}}{F} \equiv h_j(x_J) \qquad (j \in J). \qquad (4.3)$$

In practice our choice of either set *I* or *J* for checking separability will be due to the expected simplicity gained by such a choice.

**Theorem 5.** A function $F(x_1,...,x_n)$ of class $C^2(D)$ is separable in $x_I$ if and only if

$$FF_{,ij} - F_{,i}F_{,j} = 0 \qquad (i \in I, j \in J). \qquad (4.4)$$

Proof. If *F* is separable in $x_I$, then for each $(i \in I)$ $F_{,i}/F$ depends only on $x_I$, and hence $(F_{,i}/F)_{,j} = 0$ for all $(i \in I, j \in J)$. These conditions are obviously equivalent to (4.4). Conversely, if the conditions (4.4) hold, then $(F_{,i}/F)_{,j} = 0$ $(i \in I, j \in J)$, which imply that $F_{,i}/F$ are functions of $x_I$ only. By theorem (4) *F* is separable in the set $x_I$.

Let $P = \{I_1,...,I_r\}$ be a partition of $\{1,...,n\}$. Generalization of the last two theorems to treat separation with respect to the latter subsets is straightforward.

**Theorem 4g.**

The function *F*, which is of class $C^1(D)$, is separable in accordance with the partition *P* iff

$$F_{i_s}/F \equiv g_{i_s}(I_s) \qquad (i_s \in I_s, \ s = 1,...,r). \ |$$

**Theorem 5g.**

The function *F*, which is of class $C^2(D)$, is separable in accordance with the partition P iff

$$FF_{,i_s i_t} - F_{,i_s}F_{,i_t} = 0 \ (i_s \in I_s, \ i_t \in I_t; \ s<t, \ s,t=1,...,r) \ |$$

Note that theorems 1-3 corresponds to particular cases of the last two theorems. Theorem 1 corresponds to $I = \{1\}$, theorem 2 and 3 corresponds to $I_s = \{x_s\}, (s=1,...,n)$.

A good tactics to explore whether a function *F* is separable or not, and if separable, to determine its type of separability, is to form the matrix

$$M_{ij} = FF_{,ij} - F_{,i}F_{,j} \qquad (i,j = 1,...,n)$$

The function *F* will be
(i) totally separable if $M_{ij}$ is diagonal.
(ii) separable in the set of variables $x_I$ if $M_{ij}$ has a diagonal block form, when the indices $\{1,...,n\}$ are reordered so that *I* comes first.
(iii) non-separable otherwise.

## 5. Separation Using Margins.

**Definition.** Let $(x_1,...,x_n) = (x_I, x_J)$. The function $F(x_I, a_J)$ will be called the *J*-margin of *F* at $a_J$.

The following theorem shows that a function *F* is separable in the set of variables *I*, if and only if its product by its value at any point $(a_I, a_J)$ is equal to the product of its *J*-margin at $a_J$ by its *I*-margin at $a_I$.

**Theorem 6** : A function $F(x_1,....,x_n)$ is separable with respect to the set of variables $x_I$ if and only if

$$F(a_I,a_J)F(x_I,x_J) = F(x_I,a_J)F(a_I,x_J), \qquad (5.1)$$

where $(a_I,a_J)$ is a any point in the domain of $F$ at which $F(a_I,a_J) \neq 0$. Moreover, if the latter condition holds, then $F$ can be written as

$$F = G(x_I)H(x_J), \qquad (5.2)$$

where

$$G(x_I) = F(a_I,a_J)^{-1}F(x_I,a_J), \qquad H(x_j) = F(a_I,x_J). \qquad (5.3)$$

Proof. If $F$ is separable in the form (5.2), then

$$\begin{aligned} F(a_I,x_J) &= G(a_I)H(x_J) \\ F(x_I,a_J) &= G(x_I)H(a_J). \end{aligned} \qquad (5.4)$$

Multiplying the last equations side to side yields

$$F(x_I,a_J)F(a_I,x_J) = G(a_I)H(a_J)G(x_I)H(x_J) = F(a_I,a_J)F(x_I,x_J).$$

To prove the converse we assume that the conditions (5.1) hold. Dividing both sides of equations (5.1) by $F(a_I,a_J)$ yields the separable form (5.2) of $F$ in which $G$ and $H$ are given by (5.4).¦

Even $F$ is separable in a set of variables $I$, we cannot a priori specify this set. In fact we have to check the validity of the conditions (5.1) for all sets that contain one element, two elements, … and so on. The number of conditions that have to be checked is of order $\frac{1}{2}\{C_n^1 + C_n^2 +....+ C_n^{n-1}\}$, which is $2^{n-1} - 1$ conditions. Owing to the large magnitude of this number when $n$ is large, it may be better to appeal to the criterion (4.4), which requires less numbers of checking operations. For a small value of $n$ the current criterion (5.1) is indisputably the easiest.

## 6. Some Properties of Separable Functions

The following facts are almost direct consequences of the definition of separable functions.
(i) The product of two separable functions is a separable function. (The product function is defined on the intersection domain).
(ii) The reciprocal of a separable function is separable. ( If $F$ is defined on $D$, then $F$ is defined on $D - O$, where $O$ is the set of zeros of $F$. The domain $D_-$ of $1/F$ may be extended to points at which $F$ itself is not defined.

Example: $F = \sin x / \cos y$, $D = \Re - \{(2k+1)\frac{\pi}{2} : k \in Z\}$, $D_- = \Re - \{k\pi : k \in Z\}$.).

(iii) The set of separable functions $S$ defined on $D$ and have no zeros in $D$, form a group with respect to function multiplication. This follows because $S$ is closed under product and taking reciprocals, the separable function 1 acts as identity, and product is associative.

(iv) If $F$ is separable, then so are $F^n (n \in Z), |F|, |F|^\alpha (\alpha \in \Re)$.

We proceed now to consider the partial derivatives (mixed and unmixed) of a separable function $F$. It is assumed that the function $F$ is differentiable on its domain of definition at least of the order required for these expressions to be meaningful.

(v) It is easy to show that
$$\frac{\partial^k F}{\partial x_i^k} = \frac{\partial^k f_i}{\partial x_i^k} R_i = \frac{F}{f_i} \frac{\partial^k f_i}{\partial x_i^k}.$$

In particular $F_{,i} = F(f_{i,i} / f_i)$.

(vi)
$$\frac{\partial^{r+s} F}{\partial x_i^r \partial x_j^s} = \frac{F}{f_i f_j} \frac{\partial^r f_i}{\partial x_i^r} \frac{\partial^r f_j}{\partial x_j^s} \quad (i \neq j)$$

(vii) From (v) and (vi)
$$F \frac{\partial^{r+s} F}{\partial x_i^r \partial x_j^s} = \frac{\partial^r F}{\partial x_i^r} \frac{\partial^s F}{\partial x_j^s} \quad (i \neq j).$$

For $r = s = 1$, we get the conditions of separability (4.4).

(viii)
$$F^{r-1} \frac{\partial^{s_1 + \ldots + s_r} F}{\partial x_{i_1}^{s_1} \ldots \partial x_{i_r}^{s_r}} = \frac{\partial^{s_1} F}{\partial x_{i_1}^{s_1}} \cdots \cdots \frac{\partial^{s_r} F}{\partial x_{i_r}^{s_r}}$$

where $\{i_1, \ldots, i_r\}$ is a subset of $\{1, 2, \ldots, n\}$. In particular
$$F^{n-1} F_{,12\ldots n} = F_{,1} F_{,2} \ldots F_{,n}.$$

(ix) $\quad\quad\quad\quad\quad\quad (F_{,i} / F)_{,j} = (F_{,j} / F)_{,i} = 0.$

(x) Each partial derivative of a separable function of any rank (mixed and unmixed) is separable. This can be seen from the fact that differentiation with respect to a variable $x_i$ affects only the factor $f_i$. This can also be seen from (viii).

The last fact (x) constitutes a proof of the following obvious theorem.

**Theorem 7.** A function $F$ that has a certain partial derivative (of any order) which is not separable, is itself not separable. ¦

The latter theorem can be used as a tool for refuting separability.

Example: $F = x^3 y + x^2 y^2 + xy + y^2$. As we recognize $F_{,xx} = (6x + 2y)y$ as non-separable, $F$ is also non-separable.

It evident that the separability of all partial derivatives of a function are not sufficient conditions for its separability. This best can be seen through a counter example. The function $F = x^2 + y^2$ is non-separable, and yet all its partial derivatives are separable.

If $F_{,1}$ is not separable then so are all functions of the type $F + R_1(x_2,....,x_n)$. If $F_{,12}$ is not separable then so are all functions of the form $F + x_1 R_{12} + x_2 R'_{12}$, where $R_{12}$ and $R'_{12}$ are arbitrary functions independent of $x_1$ and $x_2$.

## 7. Separability and Linear Transformations

We introduce here a new set of variables $\bar{x}_i$ defined by the general linear transformations

$$\bar{x}_i = \sum_{j=1}^{n} T_{ij} x_j + b_i \qquad (i = 1,...,n) \qquad (7.1)$$

in which $T_{ij}$ is an $(n \times n)$ matrix and $b_i$ is an n-vector. A function $F(x_1,....,x_n)$ transforms under (7.1) as an invariant

$$\bar{F}(\bar{x}_1,....,\bar{x}_n) = F(x_1(\bar{x}),....,x_n(\bar{x})) \qquad (7.2)$$

where the coordinates $\bar{x}_1,....,\bar{x}_n$ are denoted collectively by $\bar{x}$.

The question posed here is that: does separability persist under the transformation (7.1)? If not, is it possible to impose additional conditions so that a separable function $F$ remains separable? And to what extent this additional conditions narrow the class of functions fulfilling this property? The answer to the first part of this question is intuitively negative. As a counter example, we mention that the separable function $F = xy$ transforms under $(x = \bar{x} + \bar{y}, y = \bar{x} - \bar{y})$ to the non-separable form $\bar{F} = \bar{x}^2 - \bar{y}^2$. To answer the second and third parts we revert to eqs(4.4), which gives the necessary and sufficient conditions of separability, and consider the transform of its left hand-side under (7.1):

$$\bar{F}\bar{F}_{,ij} - \bar{F}_{,i}\bar{F}_{,j} = \sum_{k,l=1}^{n} T_{ik} T_{jl}(FF_{,kl} - F_{,k}F_{,l}) = \sum_{k=1}^{n} T_{ik}T_{jk}(FF_{,kk} - F_{,k}F_{,k}) \quad (i \neq j),$$

since $FF_{,kl} - F_{,kl}F_{,kl} = 0$ if $(k \neq l)$. The left hand side vanishes if and only if $FF_{,kk} - F_{,k}F_{,k} = 0$ for all $k$. Thus a function $F$ has to fulfill the equations

$$FF_{,ij} - F_{,i}F_{,j} = 0 \qquad (7.3)$$

for all $i$ and $j$, including $i = j$, in order to remain separable under the transformation (7.1). In the latter case the set of equations (7.3) form a tensor equation ($F$ is an invariant, $F_{,i}$ is a vector, and $F_{,ij}$ is a tensor of second rank) with respect to the transformation (7.1). This implies that if this equation is valid with respect to one set of variables, it is valid with respect to all set of variables that result from the original set by a general linear transformation (7.1). The vanishing of the quantities (7.3) for $i = j$ is equivalent to vanishing of $(F_{,i}/F)_{,i} = (f_{i,i}/f_{,i})$, which are generally functions of $x_i$, and need not to vanish. Therefore, a separable function remains separable under the transformation (7.1) if and only if $(f_{i,i}/f_i) = 0$ for each $i$. This is equivalent to $f_i = \exp(a_i x + c_i),$ where $a_i$ and $c_i$ are arbitrary constants. Hence a function that retains its separability under the transformation (7.1) is of the form

$$F = \exp[\sum_{i=1}^{n} a_i x_i + c].$$

where $a_i \ (i = 1,...,n)$ and $c$ are arbitrary constants. This result might have been guessed from start.

In the following two section we determine the conditions imposed on the coefficients of a polynomial $P(x_1,...,x_n)$ in order to be separable.

## 8. Polynomials' Variable Separation.

As a starting point we consider the product of two polynomial $X_n(x)$ and $Y_m(y)$ of degrees $n$ and $m$ respectively:

$$X_n(x) Y_n(y) = \sum_{i=0}^{n} a_i x^i \cdot \sum_{j=0}^{m} b_j y^j = \sum_{i=0}^{n}\sum_{j=0}^{m} a_i b_j x^i y^j \equiv \sum_{i=1}^{n}\sum_{j=1}^{m} c_{ij} x^i y^j \quad (8.1)$$

The degree of the leading term $a_n b_m x^n y^m$, namely $n + m$, is called the degree of the polynomial. The other terms have degrees that descend from $(n + m - 1)$ to 0. The total number of terms in the product is $(n+1).(m+1)$ which gives rise to the same number of coefficients $c_{ij} = a_i b_j$. These coefficients however, are determined by the $n + m + 2$ values $a_i$ and $b_j$. It is clear therefore that severe limitations are imposed on the coefficients of a polynomial $P(x,y)$ so that it admits variable separation. Note that the condition $c_{nm} \neq 0$, which guarantees that the polynomial is of degree $n + m$, discards some forms of polynomials as utterly non-separable. To single out such types of anomalous polynomials, we take the term in which $i + j$ is maximum, say $m + n$. If there exists any other term $x^k y^l$ in which $k + l = n + m$, or $k > n$, or $l > m$, then the polynomial $P(x,y)$ is non-separable. For example, the polynomials $P = x^2 y^4 + x^3 y^3$, $Q = x^3 y^3 + xy^4$ are certainly not separable. A polynomial $P(x,y)$ of degree $n + m$ is denoted by $P_{nm}(x,y)$. The factors $X_n(x)$ and $Y_m(y)$ appearing in the separation of $P_{nm}(x,y)$ are unique up to constant multipliers of the form $\alpha$ and $1/\alpha$. With no loss of generality, we shall take the coefficients of the leading terms in the given polynomial and in its factors equal to 1. Thus we set

$$c_{nm} = 1, \ a_n = 1, \ b_m = 1. \quad (8.2)$$

This choice is possible since $c_{nm} = a_n b_m$. Now, the polynomial $P_{nm}(x,y)$ is separable if and only if there exists polynomials $X_n(x)$ and $Y_m(y)$ such that the following equality

$$\sum_{i=1}^{n}\sum_{j=1}^{m} c_{ij} x^i y^j = \sum_{i=1}^{n}\sum_{j=0}^{m} a_i b_j x^i y^j. \quad (8.3)$$

holds identically in $x$ and $y$. This identity holds if and if

$$c_{ij} = a_i b_j \ (i = 0,1,...,n; j = 0,1,...,m). \quad (8.4)$$

Equations (8.4) thus constitute necessary and sufficient conditions for the separability of $P_{nm}(x,y)$. The latter conditions together with the conventions (8.2) determine, if $P_{nm}$ was separable, the coefficients of the factor polynomials $X_n(x)$ and $Y_m(y)$, as

$$c_{im} = a_i b_m = a_i, \quad c_{nj} = a_n b_j = b_j. \tag{8.5}$$

Having determined the coefficients of the factors polynomials by (8.5), it suffices to fulfill the remaining relations of (8.4), namely those corresponding to $(i = 0,1,...,n-1; j = 0,1,..,m-1)$. Substituting for $a_i$ and $b_j$ from (8.5) into (8.4) yields the conditions of separability given by

$$c_{ij} = c_{im} c_{nj} \quad (i = 0,...,n-1; j = 0,...,m-1) \tag{8.6}$$

Due to the choice (8.2) one can assume that the latter conditions are already valid for $i = n$ and $j = m$. We sum up the results we have obtained above in the following theorem.

**Theorem 8.** The polynomial $P_{nm}(x, y)$ is separable if and only if its coefficients $c_{ij}$ satisfy the matrix equation

$$\begin{pmatrix} c_{n-1\,m-1} & \cdots & c_{n-1\,1} & c_{n-1\,0} \\ \vdots & \vdots & \cdots & \vdots \\ c_{1\,m-1} & \cdots & c_{11} & c_{10} \\ c_{0\,m-1} & \cdots & c_{01} & c_{00} \end{pmatrix} = \begin{pmatrix} c_{n-1\,m} \\ \vdots \\ c_{1\,m} \\ c_{0\,m} \end{pmatrix} \begin{pmatrix} c_{n\,m-1} & \cdots & c_{n1} & c_{n0} \end{pmatrix} \tag{8.7}$$

(ii) When $P_{n,m}(x, y)$ is separable, its factors polynomials are given by

$$X_n(x) = \sum_{i=0}^{n} c_{im} x^i, \quad Y_m(y) = \sum_{j=0}^{m} c_{nj} y^j. \tag{8.8}$$

Proof. (i) follows on noting that (8.7) is the matrix form of the equalities (8.6), which are themselves necessary and sufficient conditions for the separability of $P_{nm}(x, y)$.

(ii) follows from equations (8.5), which determine the coefficients of the polynomial's factors.

**Theorem 9.** (i) A polynomial $P_{nm}(x, y)$ is separable if and only if the product of the first column and first row of its coefficient matrix

$$C = \begin{pmatrix} c_{nm} = 1 & c_{n\,m-1} & \cdots & c_{n0} \\ c_{n-1\,m} & c_{n-1\,m-1} & \cdots & c_{n-1\,0} \\ \vdots & \vdots & \cdots & \vdots \\ c_{0m} & c_{0\,m-1} & \cdots & c_{00} \end{pmatrix}, \tag{8.9}$$

is equal to C itself.

(ii) If $P_{nm}(x, y)$ is separable then the coefficients of its factors $X_n(x)$ and $Y_m(y)$ are given by the first column and first row of C respectively.

Proof. (i) follows from (8.7) on noting that adjoining 1 as the first element in the column and row in the right hand-side results in enlarging the matrix in the left hand-side in (8.7) by an the additional column and row appearing in C.

(ii) follows from (8.5).

Denoting the $i$-th row and the $j$-th column of the matrix C by $R_i$ and $C_j$ respectively, there follows from (8.6) and (8.8) that each of the following relations is necessary and sufficient for $P_{nm}(x, y)$ to be separable:

(i) $$C = \begin{pmatrix} C_1 & c_{n\,m-1} C_1 & \cdots & c_{n0} C_1 \end{pmatrix}.$$

(ii)
$$C = \begin{pmatrix} R_1 \\ c_{n-1\,m} R_1 \\ \vdots \\ c_{0m} R_1 \end{pmatrix}.$$

Thus we have:
**Theorem 10.** Each of the following statements is necessary and sufficient for the polynomial $P_{nm}(x,y)$ to be separable:
(i) Each column of the coefficients' matrix $C$ is a multiple of its first column by the first element of the column concerned.
(ii) Each row of the coefficients' matrix is a multiple of the first element of this row by the first row of $C$.

An algorithm for checking whether $P_{nm}(x,y)$ is separable, and if it was, to carry out variable separation is the following:
(i) Write down the coefficients $c_{ij}$ as the $(n+1) \times (m+1)$ matrix (8.9),
(ii) Form matrix product $D_{ij}$ of the first column and first row of the matrix $C$.
(iii) If $D_{ij} \neq C_{ij}$, the separation is impossible.
(iv) If $D_{ij} = C_{ij}$, form the factors $X_n(x)$ and $Y_m(y)$ by (8.8).

**Example.** Consider the polynomial
$$P_{43}(x,y) = x^4 y^3 + 2x^4 y^2 - x^4 y + 3x^4 - 3x^3 y^3 - 6x^3 y^2 + 3x^2 y - 9x^3$$
$$+ 2xy^3 + 4xy^2 - 2xy + 6x + 7y^3 + 14y^2 - 7y + 21.$$
Since the coefficients matrix is equal to the product of its first column by its first row:
$$\begin{pmatrix} 1 & 2 & -1 & 3 \\ -3 & -6 & 3 & -9 \\ 5 & 10 & -5 & 15 \\ 2 & 4 & -2 & 6 \\ 7 & 14 & -7 & 21 \end{pmatrix} = \begin{pmatrix} 1 \\ -3 \\ 5 \\ 2 \\ 7 \end{pmatrix} (1 \quad 2 \quad -1 \quad 3),$$
the given polynomial is separable, and
$$P_{43}(x,y) = (x^4 - 3x^3 + 5x^2 + 2x + 7)(y^3 + 2y^2 - y + 3).$$
Note that the coefficients of $X_4(x)$ are the elements of the first column and those of $Y_3(y)$ are the elements of the first row in the coefficients matrix.

Finally it is noted that one may make use of margins to set up a computational scheme in which $P_{nm}(x,y)$ is separated after verifying the possibility of its separation. Indeed, and if $c_{00} \neq 0$,
$$P_{nm}(x,0) = \sum_{i=1}^{n} c_{i0} x^i, \quad P_{nm}(0,y) = \sum_{j=1}^{m} c_{0j} y^j. \qquad (8.10)$$
The factors $X_n(x)$ and $Y_m(y)$ are proportional to the latter corresponding margins. Since the coefficients of the leading terms in these factors are equal to 1, it is evident that these two margins must be divided by $c_{n0}$ and $c_{0m}$ respectively. The latter numbers are not zero if and only if $c_{00} \neq 0$, since $c_{00} = c_{n0} c_{0m}$. Thus we state on:

**Theorem 11**. If $c_{00} \neq 0$, the factors of the separable polynomial $P_{nm}(x,y)$ are

$$X_n(x) = c_{n0}^{-1} P_{nm}(x,0), \quad Y_m(y) = c_{0m}^{-1} P_{nm}(0,y) \quad (8.11)$$

It is evident that any computational scheme based on the last relations must, either verify first that $P_{nm}(x,y)$ is separable, and hence go through the necessary requirements of theorems 9 or 10, or alternatively, perform the product of the factors (8.11), to ensure that this product is identical to $P_{nm}(x,y)$.

## 9. Polynomials in n Variables

In this section we discuss total separation of polynomials of degree $N_1 + .... + N_n$ in $n$ variables

$$P_{N_1....N_n}(x_1,....,x_n) = \sum_{i_1=0}^{N_1} .... \sum_{i_n=0}^{N_n} c_{i_1.....i_n} x_1^{i_1} ..... x_n^{i_n}, \quad (9.1)$$

where $c_{N_1....N_2} \neq 0$. The factors of $P_{N_1.....N_n}$, when it is totally separable, are denoted by $X_1(x_1),....,X_n(x_n)$, where

$$X_r(x_r) = \sum_{i_r=0}^{N_r} a_{i_r} x_r^{i_r} \quad (r = 1,...,n) \quad (9.2)$$

It is evident that the coefficients of $P_{N_1....N_n}$ must relate to the coefficients of its factors by the equations

$$c_{i_1 i_2.....i_n} = a_{i_1} a_{i_2} ..... a_{i_n} \quad (i_r = 0,1,...,N_r; r = 1,....,n), \quad (9.3)$$

which constitute necessary and sufficient conditions for its total separability. We choose to take the coefficients of the leading terms in $P$ and its factors equal to 1:

$$c_{N_1....N_n} = 1, \quad a_{N_1} = ..... = a_{N_n} = 1. \quad (9.4)$$

From (9.3) and (9.4) the coefficients of the factors are given by

$$a_{i_1} = c_{i_1 N_2....N_n}, \quad a_{i_2} = c_{N_1 i_2 N_3....N_n}, \quad ........, \quad a_{i_N} = c_{N_1 N_2.....N_{n-1} i_n} \quad (9.5)$$

The necessary and sufficient conditions for separability, according to (9.3) and (9.5), are reduced to

$$c_{i_1....i_n} = c_{i_1 N_2....N_n} \cdot c_{N_1 i_2 N_3....N_n} \cdot .... \cdot c_{N_1.....N_{n-1} i_n} \quad (9.6)$$

We summarize the latter results by the following theorem:

**Theorem 12**. Consider the polynomial $P_{N_1....N_n}(x_1,...,x_n)$.
(i) $P$ is totally separable if and only if its coefficients fulfill the relations (9.6)
(ii) If $P$ is separable then the separation form is

$$P = \sum_{i=0}^{N_1} c_{i_1 N_2....N_n} x_1^{i_1} \cdot \sum_{i_2=0}^{N_2} c_{N_1 i_2 N_3....N_n} x_2^{i_2} \cdot ..... \sum_{i_n=0}^{N_n} c_{N_1....N_{n-1} i_n} x_n^{i_n} \quad (9.7)$$

It is noted that the relations (9.5) and (9.6), and no matter how large $n$ is, may be programmed easily, giving rise to a computational scheme that checks and carries out separation. For hand manipulations the following obvious theorem, which is given for three variables, proves useful.

**Theorem 13**. The polynomial $P_{IJK}(x,y,z)$ is totally separable iff

(i) The matrix $[c_{ijK}]$ is equal to the product of its first column and first row, and

(ii) The matrices $[c_{ijk}], (k = K-1,....,1,0)$ are multiples of $[c_{ijK}]$ by $c_{IJk}$.

Generalization of the latter theorem to n variables should have become clear. In order that a polynomial $P_{N_1...N_n}(x_1,...,x_n)$ is separable the following must hold:

- The product of the column $[c_{i_1 N_2...N_n}]$ by the row $[c_{N_1 i_2 N_3....N_n}]$ must equal to the matrix $[c_{i_1 i_2 N_3....N_n}]$.

- For $r = 3,....,n; i_r = 0,1,...,N_r,$ the matrices $[c_{i_1 i_2 N_3...i_r...N_n}]$ should be multiples of $[c_{i_1 i_2 N_3...N_n}]$ by $c_{N_1 N_2...i_r....N_n}$.

If the previous conditions hold, $P$ is factorized as given by (9.7).

Manipulations based on the last theorem become more clumsy as $n$ increases, and it is recommended to use the computational scheme based on (9.6).

**Example.**
$$P_{234}(x,y,z) = x^2 y^3 z^4 + 2x^2 y^3 z + x^2 yz^4 + 2x^2 yz + 2xy^3 z^4 + 4xy^3 z$$
$$+ 2xyz^4 + 4xyz + 3y^3 z^4 + 6y^3 z + 3yz^4 + 6yz \qquad (9.8)$$

According to theorem 12, $P_{234}(x,y,z)$ is separable iff
$$c_{ijk} = c_{i34} c_{2j4} c_{23k}$$

Utilizing theorem (13), we fix $k$ temporarily at its maximum value 4, make the substitution $c_{234} = 1$ in the last equations, then check that $c_{ij4} = c_{i34} c_{2j4}$ benefiting of theorem 9. Indeed

$$[c_{ij4}] = \begin{pmatrix} 1 & 0 & 1 & 0 \\ 2 & 0 & 2 & 0 \\ 3 & 0 & 3 & 0 \end{pmatrix} = \begin{pmatrix} 1 \\ 2 \\ 3 \end{pmatrix} (1 \ 0 \ 1 \ 0). \qquad (9.9)$$

To carry on the checking process, it suffices to show that the coefficients matrices corresponding to the remaining values of $k$ must appear as multiples of the matrix $[c_{ij4}]$ by $c_{23k}$. Since $c_{233} = c_{232} = c_{230} = 0$, there must appear in $P$ the term in z only, with coefficients given by the matrix $c_{231}[c_{ij4}] = 2[c_{ij4}]$. A glance at the given polynomial shows that this is the case, and hence $P$ is separable. Since only two coefficients $c_{23k}$, do not vanish, with that corresponding to $k = 1$ being equal to 2 (of course $c_{234}$ can not vanish, and it is certainly equal to unity) we have by (9.9)
$$P_{234}(x,y,z) = (x^2 + 2x + 3)(y^3 + y)(z^4 + 2z).$$

**References.**